\documentclass[titlepage,twoside,12pt]{article}
\usepackage{amssymb}
\usepackage{amsfonts}
\begin{document}

{\Large \bf PIIPTI, or the \\ \\ Principle of Increasing Irrelevance of \\ \\ Preference Type
Information} \\ \\

{\it Elem\'{e}r E Rosinger \\ Department of Mathematics \\ University of Pretoria \\ Pretoria, 0002 South Africa \\
e-mail : eerosinger@hotmail.com} \\ \\

{\bf Abstract} \\

It is shown that in the case of a {\it single decision maker} who optimizes several possibly
conflicting objectives, the amount of information available in {\it preference} relations
among pairs of possible decisions, when compared with all other possible information, is
tending to zero {\it exponentially} with the number of those different objectives.
Consequently, in the case of a larger number of conflicting objectives, the only way to obtain
a satisfactory amount of information is by the use of {\it non-preference} type relations
among possible decisions. \\ \\

{\bf 1. Introduction} \\ \\

There are three wider areas of decision making, each known to be subject to certain deep
structural {\it limitations}. \\

In {\it games}, two or more players make decisions, each pursuing his or her own best interest.
As is known, Binmore [1-3], the complexities involved are not seldom such as to lead to {\it
algorithmically unsolvable} situations. In this regard it is worth mentioning that, during the
late 1940s and early 1950s, when game theory had known a massive interest and development,
there was not much awareness about the possibility of the presence of the type of deep
difficulties which would more than three decades later be pointed out by Binmore. \\

In {\it social choice} the deeper structural difficulties came to attention relatively early
with Arrow's celebrated {\it impossibility}. Here, the issue is the appropriate aggregation of
a number of individual preferences. And as it turns out, this in general is not possible,
unless there is a "dictator". \\

Social choice can be seen as a {\it partial} version of the situation of a {\it single
decision maker}. Indeed, the aggregation of the set of individual choices amounts to a single
decision. On the other hand, even if such a decision is made by a single decision maker, he or
she is not supposed to be partial in any way with respect to any of the individual preferences
which are aggregated. And yet, Arrow's impossibility rules the realms of social choice. \\

The case of a {\it single decision maker} facing several conflicting objectives cannot -
according to one line of argument - be but {\it more difficult}, since he or she is not barred
from having specific preferences and expressing them in his or her decision. Consequently, the
possible {\it limitations} to be faced by a single decision maker may quite likely be {\it
more severe} than those which lead to Arrow's impossibility in social choice. \\
According to another line of argument, however, in the case of one single decision maker, the
fact that he or she faces all alone his or her own conflicting objectives gives an easy and
natural opportunity for certain cooperative type approaches. After all, cooperation can
involve bargaining, and in the case of one single decision maker with several conflicting
objectives, he himself, or she herself may end up as if bargaining with himself or herself. \\
In fact, certain forms of cooperation find a most appropriate context precisely within the
thinking of a single decision maker who exhibits a rational behaviour. \\

As it happens nevertheless, there seems to be little awareness in the literature about the
above mentioned issues relating to single decision makers who simultaneously face several
conflicting objectives. \\

Here we present, as one of the two main difficulties facing a single decision maker, what has
been named the Principle of Increasing Irrelevance of Preference Type Information, or in short,
PIIPTI, see Rosinger [4,5]. \\

The other difficulty, mentioned in short in Conclusions, is presented in some detail in
Rosinger [1-5]. \\ \\

{\bf 2. A Single Decision Maker with Multiple Conflicting \\ \hspace*{0.4cm} Objectives} \\ \\

In order to illustrate in detail what is involved, let us consider the following large and
practically important class of decision making situations, when the {\it single decision
maker}~ SDM has to deal with $n ~\geq~ 2$ typically conflicting objectives given by the
utility functions, see von Neumann \& Morgenstern, Luce \& Raiffa,

\bigskip
(2.1) \quad $ f_1, ~.~.~.~ , f_n ~:~ A ~~\longrightarrow~~ \mathbb{R} $

\medskip
and his or her aim is to maximize all of them, taking into account that most often such a
thing is not possible simultaneously, due to the conflicts involved. \\

The {\it main difficulty} of this situation is that the SDM is {\it not} supposed to have
available under any form whatsoever an overall utility function

\bigskip
(2.1$^*$) \quad $ f ~:~ A ~~\longrightarrow~~ \mathbb{R} $

\medskip
which would hopefully synthesize his or her position with respect to those $n \geq 2$
conflicting objectives in (2.1) taken simultaneously in their totality. \\

Here, as before, the set $A$ describes the available choices, namely, those which the SDM
has, and this set $A$ may as well be an infinite set, for instance, some open or closed
bounded domain in a finite dimensional Euclidean space. \\
Clearly, the functions $f_i$ in (2.1) can be seen as utility functions, and as such, they
generate preference relations on the set of choices $A$. Namely, the preference relation
$\leq_i$ corresponding to the utility function $f_i$ is defined for $a, b \in A$, by

\bigskip
(2.2) \quad $ a \leq_i b ~~~\Leftrightarrow~~~ f_i ( a ) \leq f_i ( b ) $

\medskip
In this way, the problem in (2.1) can be reduced to a choice, according to the {\it natural}
partial order (2.4) below, of a point in the set of all possible {\it decision outcomes}

\bigskip
(2.3) \quad $ B ~=~ \{~ ( f_1 ( a ), ~.~.~.~ , f_n ( a ) ) ~|~ a \in A ~\} ~\subseteq~
                                                                               \mathbb{R}^n $

\medskip
that is, the set of n-tuples of outcomes $( f_1 ( a ), ~.~.~.~ , f_n ( a ) ) \in \mathbb{R}^n$
which correspond to various choices $a \in A$ which the SDM can make. \\

Needless to say, the situation described by (2.1) is not the most general one, since it is
possible to encounter cases when the objectives are not given by utility functions, or simply,
are not even quantifiable. However, the model in (2.1) can nevertheless offer an edifying
enough situation, in order to be able to obtain relevant insights into the nature and extent
of the complexities and difficulties which a SDM can face. Furthermore, it can also lead to
general enough solution methods, including ways to choose solution concepts, see Rosinger
[1-5]. \\

Next we give {\it three} different arguments supporting PIIPTI. The first and the third ones
are of a {\it geometrical} nature related to finite dimensional Euclidean spaces. The second
argument is of a simple {\it probabilistic-combinatorial} kind. Here we should mention that,
while the first geometric argument is rather simple and obvious, the other geometric argument,
although quite elementary, appears however to be less well known, although it has important
connections with Physics. \\

{\bf A First Argument.} We start with a very simple geometric fact about finite dimensional
Euclidean spaces which can give a good insight into the more involved result in (2.8). On the
n-dimensional Euclidean space $\mathbb{R}^n$, with $n \geq 1$, we consider the {\it natural}
partial order relation $\leq$ defined for elements $x = ( x_1, ~.~.~.~ , x_n ),~ y =
( y_1, ~.~.~.~ ,  y_n ) \in \mathbb{R}^n$, according to

\bigskip
(2.4) \quad $ x ~\leq~ y ~~~\Leftrightarrow~~ x_i ~\leq~ y_i,~~\mbox{with}~ 1 \leq i \leq n $

\medskip
Let us denote by

\bigskip
(2.5) \quad $ P_n ~=~ \{~ x  \in \mathbb{R}^n  ~~|~~ x ~\geq~ 0 ~\} $

\medskip
the set of nonnegative elements in $\mathbb{R}^n$, corresponding to the partial order
$\leq$. \\

Then we can note that, for $n = 1$, the set $P_1$ is half of the space $\mathbb{R}^1 =
\mathbb{R}$. Further, for $n = 2$, the set $P_2$ is a quarter of the space $\mathbb{R}^2$. And
in general, for $n \geq 1$, the set $P_n$ is $1 ~/~ 2^n$ of the space $\mathbb{R}^n$. \\

It follows that in an n-dimensional Euclidean space $\mathbb{R}^n$, if one is given an
arbitrary element $x  \in \mathbb{R}^n$, then the probability for this element $x$ to be
nonnegative is $1 ~/~ 2^n$, thus it tends {\it exponentially} to zero with $n$. \\

Consequently, the same happens with the probability that two arbitrary elements $x ,~ y \in
\mathbb{R}^n$ are in the relationship $x \leq y$. Indeed, the relationship $x \leq y$ is
obviously equivalent with $0 \leq y - x$. \\

This means that one can expect a similar phenomenon to happen when trying to compare points in
the set of outcomes in (2.3), which correspond to the $n$ conflicting objectives in (2.1). \\

That very simple geometric fact is, actually, at the root of PIIPTI. \\

{\bf A Second Argument.} Let us now assume for the sake of technical simplicity that in (2.1)
we have a finite set of choices, namely

\bigskip
(2.6) \quad $ A ~=~ \{~ a_1, ~.~.~.~ , a_m ~\},~~~ m ~\geq~ 2 $

\medskip
A natural {\it single} preference relation on $A$ corresponding to (2.1), and which may try to
synthesize the respective $n$ conflicting objectives, should of course be given by a subset

\bigskip
(2.7) \quad $ S ~\subseteq~ A \times A $

\medskip
Here, for any $a,~ a^\prime \in A$, the SDM will prefer $a^\prime$ to $a$, in which case we
write $a ~\leq~ a^\prime$, or equivalently, $( a, a^\prime ) \in S$, if and only if one has
for each objective function $f_i$, with $1 \leq i \leq n$, either that $f_i ( a^\prime ) -
f_i ( a ) ~>~ 0$ and it is not negligible, or $| f_i ( a^\prime ) - f_i ( a ) |$ is
negligible. \\

Let us therefore see more precisely how much information one single preference relation $S$
can carry, when the number $n$ of conflicting objectives in (2.1) becomes large, and even if
only moderately so. This can be done quite easily by noting that in {\it typical} situations,
we can have the relation

\bigskip
(2.8) \quad $ \mbox{car}~ S ~/~ \mbox{car}~ ( A \times A ) ~=~ O ( 1 ~/~ 2^n ) $

\medskip
where for a finite set $E$ we denoted by "$\mbox{car}~ E$" the number of its elements. \\

The proof of (2.8) goes as follows, by using a combinatorial-probabilistic type argument. Let
us take any injective function $g : A \longrightarrow \mathbb{R}$, and denote by

\bigskip
(2.9) \quad $ S_g ~=~ \{~ ( a, a^\prime ) \in A \times A ~~|~~ g ( a ) ~\leq~ g ( a^\prime )
                                 ~\} $

\medskip
which is its corresponding preference relation on $A$. Then obviously

\bigskip
(2.10) \quad $ \mbox{car}~ S_g ~=~ m ( m + 1 ) ~/~ 2 $

\medskip
Now given any subset $S \subseteq A \times A$, let us denote by $P ( S )$ the probability that for an arbitrary pair $( a,
a^\prime ) \in A \times A$, we have $( a, a^\prime ) \in S$. Then clearly

\bigskip
(2.11) \quad $ P ( S_g ) ~=~ ( 1 + 1 ~/~ m ) ~/~ 2 $

\medskip
Let us assume about the objective functions in (2.1) the following

\bigskip
(2.12) \quad $ f_1, ~.~.~.~ , f_n ~~\mbox{are injective} $

\medskip
and furthermore, that their corresponding sets of preferences

\bigskip
(2.13) \quad $ S_{f_1}, ~.~.~.~ , S_{f_n} ~~\mbox{are probabilisitically independent} $

\medskip
Then we obtain, see (2.11)

\bigskip
(2.14) \quad $ P  ( S_{f_1} \bigcap ~.~.~.~ \bigcap S_{f_n} ) ~=~ P ( S_{f_1} ) ~.~.~.~
                                               P ( S_{f_n} ) ~=~ ( 1 + 1 ~/~ m )^n ~/~ 2^n $

\medskip
And now (2.8) follows, provided that $n$ in (2.1) and $m$ in (2.6) are such that

\bigskip
(2.15) \quad $ ( 1 + 1 ~/~ m )^n ~=~ O ( 1 ) $

\medskip
which happens in many practical situations. \\

As for the independence condition (2.13), let us note the following. Let us assume that the
objectives $f_1$ and $f_2$ are such that for $a,~ a^\prime \in A$ we have

\bigskip
(2.16) \quad $ f_1 ( a ) ~<~ f_1 ( a^\prime ) ~~~\Leftrightarrow~~~ f_2 ( a ) ~<~
                                                                   f_2 ( a^\prime ) $

\medskip
then obviously $S_{f_1} = S_{f_2}$, hence (2.14) may fail. But clearly, (2.16) means that
$S_{f_1}$ and $S_{f_2}$ are {\it not} independent. In the opposite case, when

\bigskip
(2.17) \quad $ f_1 ( a ) ~<~ f_1 ( a^\prime ) ~~~\Leftrightarrow~~~ f_2 ( a^\prime ) ~<~
                                                                            f_2 ( a ) $

\medskip
then obviously

\bigskip
(2.18) \quad $ S_{f_1} \bigcap S_{f_2} ~=~ \{~ ( a, a ) ~|~ a \in A ~\} $

\medskip
and (2.14) may again fail. However (2.18) once more means that $S_{f_1}$ and $S_{f_2}$ are
{\it not} independent, since they are in total conflict with one another. \\

{\bf A Third Argument.} For the sake of simplicity, let us assume that the set $B$ of outcomes
in (2.3) is of the form

\bigskip
(2.19) \quad $ B ~=~ \bigg \{~ b = ( b_1, ~.~.~.~ , b_n ) \in \mathbb{R}^n ~~ \bigg |~~
             \begin{array}{l}
                  b_1, ~.~.~.~ , b_n ~\geq~ 0 \\
                  b_1 + ~.~.~.~ + b_n ~\leq~ L
             \end{array} ~\bigg \} $

\medskip
for a certain $L > 0$. Then clearly the Pareto maximal, or in other words, the non-dominated
subset of $B$ is

\bigskip
(2.20) \quad $ B^P ~=~ \{~ b = ( b_1, ~.~.~.~ , b_n ) \in
\mathbb{R}^n ~~|~~ b_1 + ~.~.~.~ + b_n = L ~\} $

\medskip
when considered with the natural partial order (2.4) on $\mathbb{R}^n$. \\

Now for $0 < \epsilon < L$, the $\epsilon$-thin shell in $B$ corresponding to $B^P$ is given
by

\bigskip
(2.21) \quad $ B^P ( \epsilon ) ~=~ \bigg \{ b = ( b_1, ~.~.~.~ , b_n ) \in
\mathbb{R}^n ~\bigg | \begin{array}{l}
                         L - \epsilon ~\leq~  b_1 + ~.~.~.~ + b_n ~\leq~ \\
                            ~~~~~~~~~~~~~~~~~~~~~~ \leq~ L
                       \end{array} \bigg \} $

\medskip
And a standard multivariate Calculus argument gives for the volume of $B$ in (2.19) the
relation

\bigskip
(2.22) \quad $ \mbox{vol}~ B ~=~ K_n~ L^n $

\medskip
where the constant $K_n > 0$, involving the Gamma function, does only depend on $n$, but not
on $L$ as well. In this way it is easy to see that

\bigskip
(2.23) \quad $ \mbox{vol}~ B^P ( \epsilon ) ~/~ \mbox{vol}~ B ~=~ 1 - ( 1 - \epsilon ~/~ L )^n $

\medskip
This leads to a rather counter-intuitive and somewhat paradoxical property of {\it higher}
dimensional Euclidean spaces. \\
For instance, in the 20-dimensional case, a shell with a thickness of only 5\% of the radius
$L$ of a sphere will nevertheless contain at least 63\% of the total volume of that sphere. \\

In more simple and direct {\it geometric} terms the relation (2.23) means that : \\

{\bf ( VOL )} "The volume of a multidimensional solid is mostly concentrated next to its
surface." \\

The {\it relevance} of this property {\bf ( VOL )} to PIIPTI is as follows. The set (2.3), or
equivalently (2.19), of outcomes $B$ in the multiple objective decision problem (2.1) is the
one which determines the choice of the appropriate decision taken in $A$, and it does so
through the relations (2.2). And obviously, in this respect, only the Pareto maximal, or the
non-dominated subset $B^P$ of $B$, see (2.20), is of relevance. However, within this subset
$B^P$ {\it no} two different points $u \neq v$ can be in a relation $u \leq v$, see (2.4),
this being the very definition of a Pareto maximal, or non-dominated set. In this way {\it all}
pairs of different elements $u \neq v$ in $B^P$ are {\it incomparable}, thus are outside of
being included in a preference relationship. \\
And as seen in (2.23), the volume of no matter how thin a shell next to $B^P$, when compared
to that of $B$, tends to 1 exponentially with $n$ becoming large. Thus, the amount of pairs of
different elements $u \neq v$ in no matter how small a neighbourhood of $B^P$ tends to 1 with
$n$ becoming large, when compared with all the possible pairs of outcomes in $B$. And
obviously, any two different elements $u \neq v$ in a neighbourhood of $B^P$ are comparable,
that is, satisfy the relation $u \leq v$, see (2.4), only if the are very near to one another,
thus they {\it cannot} express any kind of a more relevant preference. Otherwise, if two
different elements $u \neq v$ in a neighbourhood of $B^P$ are not near to one another, then
they must be {\it incomparable}, given the fact that $B^P$ is a Pareto maximal, or
non-dominated set. \\

It may be instructive to note that relation (2.23) also has a {\it physical} interpretation,
as it explains the phenomenon of {\it temperature}, see Manin. Indeed, let us assume that a
certain simple gas has $n$ atoms of unit mass. Then their {\it kinetic energy} is given by

\bigskip
(2.24) \quad $ E ~=~ \Sigma_{1 \leq i \leq n}~ v_i^2 ~/~ 2 $

\medskip
where $v_i$, with $1 \leq i \leq n$, are the velocities of the respective atoms. Therefore,
for a given value of the kinetic energy $E$, the state of the gas is described by the vector
of $n$ velocities, namely

\bigskip
(2.25) \quad $ v ~=~ ( v_1, ~.~.~.~ , v_n ) \in S_n ( \sqrt ( 2 E ) ) $

\medskip
where for $L > 0$, we denoted by

\bigskip
(2.26) \quad $ S_n ( L ) ~=~ \{~ x = ( x_1, ~.~.~.~ , x_n ) \in
\mathbb{R}^n ~~|~~ x_1^2 + ~.~.~.~ + x_n^2 = L^2 ~\} $

\medskip
the ( n - 1 )-dimensional surface of the n-dimensional ball with radius $L$ in $\mathbb{R}^n$. \\

Now we can recall that in view of the Avogadro number, under normal conditions for a usual
macroscopic volume of gas, one can have

\bigskip
(2.27) \quad $ n ~>~ 10^{20} $

\medskip
Therefore, the above property {\bf ( VOL )} which follows from (2.23) is very much manifest.
Let us then assume that a small thermometer with a thermal energy $e$ negligible compared to
$E$ is placed in the gas. Then the state (2.25) of the gas will change to a new state

\bigskip
(2.28) \quad $ v ~=~ ( v_1, ~.~.~.~ , v_n ) \in S_n ( \sqrt ( 2 E^\prime ) ) $

\medskip
However, in view of property {\bf ( VOL )}, it will follow with a high probability that

\bigskip
(2.29) \quad $ E^\prime ~\approx~ E $

\medskip
And it is precisely this {\it stability} or {\it rigidity} property (2.29) which leads to the
phenomenon of {\it temperature} as a macroscopically observable quantity. \\ \\

{\bf 3. Conclusions} \\ \\

Situations involving the actions of conscious rational agents are approached in three
mathematical theories, namely, the theory of games, the theory of social choice, and decision
theory. In games, there are two or more such conscious and rational agents, called players,
who are interacting according to the given rules. And except for that, they are free and
independent, and there is no overall authority who could influence in any way the players. In
social choice, again, there are two or more conscious and rational agents with their  given
individual preferences. Here however, the issue is to find a mutually acceptable aggregation
of those preferences. And such an aggregation is seen as being done by an outsider. Finally,
decision theory can be seen as a two person game, in which one of the players is a conscious
rational agent, while the other is Nature. \\

In the case of one single decision maker, who in decision theory is seen as a conscious
rational player, playing alone against Nature, one may seem at first two have a situation
which enjoys all the advantages that are missing both in games, and in social choice. Indeed,
it may at first appear that such a single decision maker does not have to put up with one or
more other autonomous players. And also, as the single player, he or she can automatically be
seen as a dictator as well, since there is no other conscious agent out there to protest,
least of all what is called Nature in such a context. \\
It would, therefore, appear that in decision theory one has it rather easy. \\

And yet, in the typical practical situations when the single decision maker is facing multiple
and conflicting objectives, all the mentioned seeming advantages are instantly cancelled.
Instead, the single decision maker can easily end up feeling as if two or more autonomous
agents have moved inside of him or her, and now he or she has to turn into a dictator who, in
fact, ends up fighting himself or herself. \\

In this regard, two facts come to the fore from the beginning in typical situations with a
single decision maker facing multiple conflicting objectives, see Rosinger [1-5] : \\

{\bf Fact 1.~} There is {\it no}, and there {\it cannot} be a unique natural canonical
candidate for the very {\it concept of solution}. And in fact, the very issue of choosing a
solution concept leads to a {\it meta-decision} problem which itself has multiple conflicting
objectives. \\

{\bf Fact 2.~} The information contained in the preference structures involved - relative to
all other possible, such as for instance, non-preference type information present in the
situation - tends {\it exponentially} to zero, as the number of conflicting objectives
increases. This phenomenon, which in fact is of a very simple {\it higher dimensional
geometric} nature, can be called the Principle of Increasing Irrelevance of Preference Type
Information, or in short PIIPTI. \\ \\

\end{document}